\documentclass[12pt, letterpaper]{article}
\usepackage{lineno}
\usepackage{hyperref}
\modulolinenumbers[5]

\usepackage{subfigure}
\usepackage{tikz,gincltex}
\usepackage{arydshln}
\usepackage{nccmath}
\usepackage{mathtools}
\usepackage{algorithm}
\usepackage{algorithmic}
\usepackage{amssymb, amsmath}
\usepackage{natbib}
\usepackage{ stmaryrd }
\usepackage{wrapfig}

\newenvironment{myproof}{%
    \par\medskip\noindent%
    \textit{Proof.}\enspace\ignorespaces%
}{%
    \hfill$\square$\par\medskip%
}

\usepackage{breakcites}

\usepackage[font=small,labelfont=bf]{caption}

\usepackage{pgfplots,tikz-3dplot}
\pgfplotsset{compat=newest}

\usepackage{algorithmic}
\usepackage{pgfplots}
\usetikzlibrary{positioning,chains,fit,shapes,calc}
\usetikzlibrary{backgrounds}
\usetikzlibrary{arrows.meta}
\usetikzlibrary{shapes,snakes}

\newtheorem{theorem}{Theorem}
\newtheorem{lemma}{Lemma}
\newtheorem{proposition}{Proposition}

\newtheorem{corollary}{Corollary}
\newtheorem{remark}{Remark}
   
\newtheorem{assumption}{Assumption}

\definecolor{ao}{rgb}{0.55, 0.71, 0.0}
\definecolor{bleudefrance}{rgb}{0.19, 0.55, 0.91}
\definecolor{dimgray}{rgb}{0.41, 0.41, 0.41}    
\definecolor{mediumorchid}{rgb}{0.73, 0.33, 0.83}
\definecolor{mediumtealblue}{rgb}{0.0, 0.33, 0.71}
\definecolor{harvestgold}{rgb}{0.85, 0.57, 0.0}
\definecolor{blue(pigment)}{rgb}{0.2, 0.2, 0.6}
\definecolor{forestgreen(traditional)}{rgb}{0.27, 0.35, 0.27}
\definecolor{cadmiumred}{rgb}{0.89, 0.0, 0.13}
\definecolor{orange(webcolor)}{rgb}{1.0, 0.5, 0.0}
\definecolor{tangerine}{rgb}{0.95, 0.52, 0.0}
\definecolor{c1}{rgb}{0.368417, 0.506779,0.709798}
\definecolor{c2}{rgb}{0.880722, 0.611041,0.142051}
\definecolor{c3}{rgb}{0.560181, 0.691569,0.194885}
\definecolor{c4}{rgb}{0.922526, 0.385626,0.209179}
\definecolor{c5}{rgb}{0.528488, 0.470624,0.701351}
\definecolor{c6}{rgb}{0.772079, 0.431554,0.102387}
\definecolor{c7}{rgb}{0.363898, 0.618501,0.782349}
\definecolor{c8}{rgb}{1, 0.75, 0}
\definecolor{c9}{rgb}{0.647624, 0.37816, 0.614037}
\definecolor{c10}{rgb}{0.571589, 0.586483, 0.}
\definecolor{c11}{rgb}{0.915, 0.3325, 0.2125}
\definecolor{c12}{rgb}{0.83, 0.46, 0.}
\definecolor{c13}{rgb}{0.9575, 0.545, 0.11475}
\definecolor{c14}{rgb}{1., 0.7575, 0.}
\definecolor{c15}{rgb}{0.6175, 0.715, 0.}
\definecolor{c16}{rgb}{0.15, 0.715, 0.595}
\definecolor{c17}{rgb}{0.3625, 0.545, 0.85}
\definecolor{c18}{rgb}{0.575, 0.4175, 0.85}
\definecolor{c19}{rgb}{0.677, 0.358, 0.595}
\definecolor{c20}{rgb}{0.7875, 0.358, 0.425}
\definecolor{c21}{rgb}{0.915, 0.3325, 0.2125}

\title{Average consensus with resilience and privacy guarantees without losing accuracy\footnote{This work is supported by the Swedish Research Council under the
grant 2021-06316 and by the Swedish Foundation for
Strategic Research. }}

\author{Guilherme Ramos\thanks{Instituto Superior Técnico, Universidade de Lisboa, 1049-001 Lisbon, Portugal, and Instituto de Telecomunicações, 1049-001 Lisboa, Portugal}\qquad Daniel Silvestre\thanks{School of Science and Technology, NOVA University of Lisbon (FCT/UNL), 2829-516 Caparica, Portugal, COPELABS, Lusófona University, Lisboa, Portugal, and Institute for Systems and Robotics, Instituto Superior T\'ecnico, Universidade de Lisboa, Portugal}\\André M. H. Teixeira\thanks{Division of Systems and Control, Department of Information Technology, Uppsala University, Sweden}\qquad Sérgio Pequito\thanks{Instituto Superior Técnico, Universidade de Lisboa, 1049-001 Lisbon, Portugal, and Institute for Systems and Robotics, Instituto Superior T\'ecnico, Universidade de Lisboa, Portugal}}

\begin{document}

\maketitle
          

\begin{abstract}
This paper addresses the challenge of achieving private and resilient average consensus among a group of discrete-time networked agents without compromising accuracy. 
State-of-the-art solutions to attain privacy and resilient consensus entail an explicit trade-off between the two with an implicit compromise on accuracy. In contrast, in the present work, we propose a methodology that avoids trade-offs between privacy, resilience, and accuracy. 
We design a methodology that, under certain conditions, enables non-faulty agents, i.e., agents complying with the established protocol, to reach average consensus in the presence of faulty agents, while keeping the non-faulty agents' initial states private. 
For privacy, agents strategically add noise to obscure their original state, while later withdrawing a function of it to ensure accuracy. 
Besides, and unlikely many consensus methods, our approach does not require each agent to compute the left-eigenvector of the dynamics matrix associated with the eigenvalue one. 
Moreover, the proposed framework has a polynomial time complexity relative to the number of agents and the maximum quantity of faulty agents. Finally, we illustrate our method with examples covering diverse faulty agents scenarios.
\end{abstract}


\section{Introduction}\label{sec:intro}


Distributed consensus assumes paramount significance across a multitude of application domains where coordination and concurrence among numerous entities are crucial. 
Examples include task distribution, load equilibrium, and data coherence~\cite{chiang2007layering},  
%
robotics and autonomous to facilitate a cohort of robots or agents to 
collectively arrive at decisions, enabling collaborative undertakings such as formation management, exploration, and cooperative surveillance~\cite{jadbabaie2003coordination,alessandretti2019optimization,cortes2006robust,ribeiro:2020,ribeiro:c20}. 
Additionally, in sensor networks and Internet-of-things (IoT) applications, consensus mechanisms empower effective data fusion, anomaly detection, and decision-making among dispersed sensors or IoT devices~\cite{silvestre:desync}. 
Furthermore, consensus assumes a pivotal role in the landscape of distributed optimization~\cite{tsitsiklis1986distributed,johansson2008subgradient}. 
On the whole, distributed consensus occupies a central position across a wide array of fields, fostering efficient coordination, resilience against faults, scalability, and robustness within intricate and distributed systems. 

The relevance of consensus methods in a wide range of applications created the need to study and improve properties that go beyond accuracy. 
Namely, it is paramount to ensure properties such as \textit{privacy} and \textit{resilience} to attacks in a plethora of consensus applications. 


Henceforth, \textbf{privacy} is a crucial aspect of distributed consensus, safeguarding sensitive information and maintaining confidentiality throughout the consensus process within a distributed system~\cite{such2014survey, pequito2014design, gupta2017privacy}. 
Different approaches for privacy-preserving consensus methods are compared in \cite{ramos2023designing}. 

The present work aligns with the observability-based (O-based) strategies that focus on preventing curious agents from inferring other agents' states by observing some measurements of the dynamics’ state. In particular, our work relates to state and input estimation, where our aim is to prevent the estimation of agents' initial states~\cite{hui2005observer,reed2022minimum}, i.e., to protect agents' private states from being retrieved or estimated. Observability in dynamical systems establishes necessary and sufficient conditions for designing an estimator capable of recovering agents’ initial states, which are meant to be private~\cite{pequito2014design, boutat2021observability, ramos2021distributed, ramos2023designing}. Approaches to achieving O-based privacy may involve network augmentation~\cite{ramos2023designing, ramos2021distributed, ramos2023trade}.  
\citet{ramos2023discrete} introduce a discrete-time reputation-based consensus method, along with a continuous-time variation~\cite{ramos2022resilient}, to address attack scenarios unattainable with previous techniques. This method minimizes the impact of undetected attackers. \cite{shang2023resilient} extends these works to include multi-hop communication and path-dependent heterogeneous delays. 

Hereafter, we attain privacy by scheduling the introduction of noise over a finite window of time, which will enable us to maintain both privacy and accuracy by later removing the noise in a strategical manner, following the proposed strategy to achieve privacy in~\cite{manitara2013privacy}. While our approach is inspired by the idea of adding noise to enhance privacy, it does not employ differential privacy techniques. Instead, it leverages concepts from dynamical systems and noise scheduling to ensure that private information is obscured from potential observers.

On the other hand, we need to ensure \textbf{resilience}, the ability of a distributed system to reach a consistent agreement even in the face of failures or malicious attacks on its components. 
Resilient distributed consensus algorithms are designed to ensure that the system remains reliable and functional, even in the presence of \emph{faulty} agents, i.e., when some nodes may fail, exhibit faulty behavior, or try to disrupt the consensus process \cite{oksuz2018distributed, sundaram2018distributed, dibaji2015consensus, kikuya2017fault, sundaram2016ignoring, usevitch2018resilient, saldana2017resilient, dibaji2017resilient, chen2018attack, dibaji2019resilient}.


Previous resilience consensus methods often use the \mbox{mean-subsequence-reduced} algorithm (MSR), which converges when both the network of non-faulty agents is robust against $f+1$ faulty agents and the non-faulty agents have an equal number of in-neighbors. 
The key idea behind the MSR is that each normal agent sorts the neighbors' received values by the distance to its state and discard the $f$ largest ones.

An alternative approach to achieving resilient consensus suggests incorporating the concept of reputation in order to disregard neighboring agents with low reputation scores. 
Several studies \cite{li2012robust, saude2017robust, saude2017reputation, 10.1145/3397271.3401278} have proposed reputation-based ranking systems, demonstrating their effectiveness in mitigating attacks and bribery. 

The \textbf{main contribution} of this work is an average consensus method that synergetically guarantees resilience and privacy without compromising accuracy. 
Furthermore, our proposal avoids the computation
of eigenvectors, which typically necessitates
numerical methods and introduces additional sources of
inaccuracies. In contrast to most of the proposed methods in the literature, which compute the left-eigenvector associated with the unit eigenvalue, our method preserves the true distributed nature of the proposed protocols. 
Furthermore, the previous state-of-the-art methods~\cite{ramos2021distributed,ramos2023trade} require each node to augment, \emph{a priori}, the network by creating a local virtual network, resulting in a higher computational complexity compared to our approach.

\section{Problem statement}\label{sec:prob_stat}

Consider a discrete-time consensus method with dynamics 
\begin{equation}\label{eq:lti}
x^{(k+1)}=Ax^{(k)} + \mathbb I_{[0,T]}(k)\eta^{(k+1)},
\end{equation}
where $x^{(k)}=[x_1^{(k)} \,\ldots\, x_n^{(k)}]^\intercal$ is a vector of scalar states, each one associated with an agent, $\mathbb I_{[0,T]}:\mathbb Z_0^+\to\{0,1\}$ is a indicator function, i.e., $\mathbb I_{[0,T]}(k)=1$ if and only if $k\in\mathbb Z_0^+$ and $k\leq T$, and $\eta^{(k)}\in\mathbb R^n$ is the noise such that $\eta_i^{(k)}\sim \mathcal N(0, \xi^2)$, and, at the realization, $\eta^{(k)}_i$ is known only to agent $i$. 
Additionally, there may be potentially (at most) $f$ faulty agents, i.e., that do not obey the dynamics in~\eqref{eq:lti}.  
The goal we aim to pursue is two-fold: (i) consensus among non-faulty agents, i.e., for each non-faulty agent $i$, 
\begin{equation}\label{eq:final_consensus}
\lim_{k\to\infty}x_i^{(k)} = \frac{1}{n-|\mathcal F|}\sum_{j\in[n]\setminus\mathcal F} x_j^{(0)},
\end{equation}
where $\mathcal F\subset \mathcal V = [n]$, the set of faulty agents' indices, with $[n]=\{k\,:\,(k\in\mathbb Z)\land(1\leq k\leq n)\}$, and $|\mathcal F|\leq f$, and (ii) privacy of the state of non-faulty agents, i.e., its states cannot be (exactly) recovered by other agents. In other words, we want to ensure uncertainty associated with a non-faulty agent's initial state retrieved. 

\begin{remark}
    We do not consider passive and curious agents as faulty agents. Our set of faulty agents $\mathcal F$ only accounts for agents that are not following the consensus update rule. In fact, every agent may try to recover the initial state of other agents by observing the values shared within the network, possibly  in collusion.\hfill$\blacktriangle$
\end{remark}

For technical reasons, in what follows, we consider the following assumptions.

\begin{assumption}\label{ass:-1}
    For simplicity, we assume that the network is undirected, meaning communication between agents is bidirectional.~\hfill$\star$
\end{assumption}

\begin{assumption}\label{ass:0}
    Each agent in the network is aware of the number of agents and the agents' indices, but only knows its connections to its immediate neighbors. A curious or faulty agent may know the whole network.~\hfill$\star$
\end{assumption}

Notwithstanding, we point out in Remark~\ref{rem:doubly} under which conditions our results hold also for directed networks.  



The structural connections between agents in a given network are represented by a \textit{graph}, which is an ordered pair $\mathcal{G} = (\mathcal{V}, \mathcal{E})$, where $\mathcal{V}$ is a nonempty set of nodes (i.e., agents), and $\mathcal{E} \subseteq \mathcal{V} \times \mathcal{V}$ is a set of edges. Each edge is a pair representing the accessibility relationship between nodes. That is, if $(u, v) \in \mathcal{E}$ (and so $(v, u)\in \mathcal{E}$), then node $v$ directly accesses information from node $u$, and vice-versa, and $v$ is a neighbor of $u$ and $u$ a neighbor of $v$. 

A common representation of a graph with $n$ nodes is its adjacency matrix $A \in \mathbb{R}^{n\times n}$, where $A_{u,v} = 1$ if $(u,v) \in \mathcal{E}$, and $A_{u,v} = 0$ otherwise. 
Similarly, we associate to $A\in\mathbb R^{n\times n}$ a graph representation $\mathcal G(A) = (\mathcal V,\mathcal E)$, where $\mathcal V= [n]$ and $\mathcal E=\{(i,j)\,:\,A_{ij}\neq 0\}$. 
A \textit{subgraph} or subnetwork $\mathcal{H} = (\mathcal{V}', \mathcal{E}')$ of a graph $\mathcal{G} = (\mathcal{V}, \mathcal{E})$ is a graph such that $\mathcal{V}' \subset \mathcal{V}$ and $\mathcal{E}' \subset (\mathcal{V}' \times \mathcal{V}' \cap \mathcal{E})$. We use  $\mathcal H\equiv \mathcal{G} \setminus \mathcal{A}$, with $\mathcal{A} \subset \mathcal{V}$ to denote $\mathcal{H} = (\mathcal{V} \setminus \mathcal{A}, \mathcal{E}')$, where $\mathcal{E}' = \{(u,v) \in \mathcal{E} : u,v \notin \mathcal{A}\}$.\hfill
 
We also make the following assumption.

\begin{assumption}\label{ass:1}
    There are no $\mathcal S,\mathcal S'\subset\mathcal V\setminus\mathcal F$, with $\mathcal S\neq\mathcal S'$ and $|\mathcal S'|,|\mathcal S|\leq f$, such that 
    \[
    \frac{1}{n-|\mathcal S|}\sum_{j\in[n]\setminus\mathcal S} x_j^{(0)} = \frac{1}{n-|\mathcal S'|}\sum_{j\in[n]\setminus\mathcal S'} x_j^{(0)},
    \]
    i.e., there are no two subsets such that the average of the initial values of agents in $[n]\setminus \mathcal S$ and the average of initial states of agents in $[n]\setminus \mathcal S'$ is the same.~\hfill$\star$
\end{assumption}

Assumption~\ref{ass:1} is a mild assumption since if the agents initial states are selected at random, such a scenario occurs with zero probability~\cite{ramossilvestresIJoC}. 
Also, it is unlikely that the average of any subset of initial states is the same even when the initial states are not drawn at random. 
This assumption will be used in the soundness proof of the proposed methodology -- see Theorem~\ref{th:sound}.

Next, we detail the faulty agents model adopted in this work. 



\vspace{2mm}

Regarding the faulty agents, we do the following assumption. 

\begin{assumption}\label{ass:faulty}
Given a network of agents $\mathcal G=(\mathcal V, \mathcal E)$ and a set of faulty agents $\mathcal F\subset\mathcal V$, the following yield:
\begin{description}
    \item[$(i)\,$] the network $\mathcal H=(\mathcal V\setminus\mathcal F, \mathcal E')$ is strongly connected;
    \item[$(ii)$]  for every $v\in\mathcal F$,  
$\displaystyle\lim_{k\to\infty}x_v^{(k)} = \varphi_v,$ 
with $\varphi_v\in\mathbb R$, and, for $v\neq v'$, we can have $\varphi_v\neq \varphi_{v'}$ or $\varphi_v=\varphi_{v'}$ (each faulty agent's state converges to a constant value, possibly distinct values from other faulty agent).\hfill$\star$
\end{description}
\end{assumption}


The described scenario models various types of agents, including faulty agents that are stubborn, malicious agents attempting to steer the consensus toward a particular value, and colluding agents working together to achieve a shared target value.

{
\begin{assumption}\label{ass:priv}
Given a network of agents $\mathcal G=(\mathcal V, \mathcal E)$ and a set of faulty agents $\mathcal F\subset\mathcal V$, the following holds:
 the network $\mathcal H=(\mathcal V\setminus\mathcal F, \mathcal E')$ is such that for any node $v\in (\mathcal V\setminus\mathcal F)$ we have that $\mathcal N_v\neq \mathcal V\setminus\mathcal F$.
\hfill$\star$
\end{assumption}

Assumption~\ref{ass:priv} is crucial for maintaining privacy in the consensus protocol. If an agent were directly connected to all others, it could trivially compute and share the exact average of all states, bypassing the privacy-preserving mechanisms. This would allow for the derivation of initial states, compromising the protocol's privacy guarantees.
}



\section{Average consensus with resilience and privacy guarantees without losing accuracy}\label{sec:main_gen} 

In this section, we describe several algorithms (Algorithms~\ref{alg:ds}–\ref{alg:ss}) that serve as building blocks for the proposed methodology to achieve resilient average consensus with privacy guarantees, with accuracy, as detailed in Algorithm~\ref{alg:main}, which presents an accurate average consensus method with both resilience and privacy guarantees
without compromising on either property. 

The overall framework is as follows: Algorithm~\ref{alg:ds} designs the weights for the dynamics to ensure that average consensus is achievable -- \textbf{accurate average consensus}; Algorithm~\ref{alg:gn} introduces noise into an agent’s state while keeping track of the accumulated noise that is subsequently discounted -- \textbf{private average consensus}; and, finally, Algorithm~\ref{alg:ss} plays a crucial role in the proposed method by allowing an agent to select its state from a set of possible states -- \textbf{resilient average consensus}. 

Finally, 
in Theorem~\ref{th:sound} and Theorems~\ref{th:private} and~\ref{th:privateAll}, we provide the theoretical guarantees for the accurate average consensus
method with both resilience and privacy guarantees, respectively. 
Notably, such properties do not interfere with each other. Additionally, in Proposition~\ref{prop:complex}, we provide the polynomial time complexity relative to the number of agents and the
maximum quantity of faulty agents. 
We will now provide a more detailed explanation of the idea behind each of these building blocks.\hfill 


\subsection{Design to achieve average consensus}

First, the concept behind Algorithm~\ref{alg:ds}, in Section~\ref{sub:sai}, enables each non-faulty agent to calculate the weights for the weighted averages of its neighbors’ states using only the information from its neighbors (alternative methods can be considered, as mentioned in Remark~\ref{rem:doubly}). 

To this end, we need the following notation. 
We denote the $i$th row of $A$ by $A_i$ and the $j$th column of the $i$th row of $A$ by $A_{ij}$, where $i$ and $j$ are valid indices for the dimensions of $A$. 
We denote the parts of a set (the subsets of a set) $\mathcal X$ by $\mathcal P(\mathcal X)=\{\mathcal Y\,:\,\mathcal Y\subset \mathcal X\}$. Moreover, for $k\geq 0$, we denote by $\mathcal P(\mathcal X,k)=\{\mathcal Y\,:\,\mathcal Y\in \mathcal P(\mathcal X), |\mathcal Y|=k\}$ (the subsets of $\mathcal X$ with $k$ elements), and $\mathcal P(\mathcal X, i, k) = \bigcup_{j=i}^k\mathcal P(\mathcal X, j)$ is the subsets of $\mathcal X$ with at least $i$ elements and at most $k$ elements. 
The in-degree of a node $v \in \mathcal{V}$, denoted as $d_v$, is the number of proper neighbors of $v$, given by $d_v = |\mathcal{N}_v \setminus \{v\}|$. Similarly, the out-degree of a node $v \in \mathcal{V}$, denoted as $o_v$, is the number of nodes for which $v$ is a neighbor, i.e., $o_v = |\{u : v \in \mathcal{N}_u \setminus \{u\}\}|$. In the case of a graph, the \mbox{in-degree} is equal to the \mbox{out-degree}, and we refer to either as the node degree. 

Subsequently, the weights each agent compute to its neighbors are encapsulated in an augmented matrix $\tilde{A}$, which is a \emph{block-diagonal matrix} where each block has the dynamics of the each subnetwork of agents $\mathcal S\subset\mathcal P([n],0,f)$, which we denote by $A^{\mathcal S}$, where $f$ is the maximum number of faulty agents allowed. More precisely, given the network of agents $\mathcal G=(\mathcal V,\mathcal E)$, we have that
\[
 \tilde A= \text{blockdiag}\left(A^{\mathcal S_1},\ldots,A^{\mathcal S_r}\right) \,\text{ and }\, A^{\mathcal S_i}=\widehat{A(\mathcal G\setminus\mathcal S_i)},
\]
for $i=1,\ldots,r$ and $r=|\mathcal P([n],0,f)|$, where $A(\mathcal G\setminus\mathcal S)$ is the adjacency matrix of the subgraph $\mathcal G\setminus\mathcal S$, where the rows and columns with indices in $\mathcal S$ are all set to zero to simplify the exposition and to avoid renumbering agents ids. Additionally, $\,\widehat{\,\quad\,}:\mathbb R^{n\times n}\to\mathbb R^{n\times n}$ is a function that returns a doubly stochastic matrix with the same zero entries as $A(\mathcal G\setminus\mathcal S)$. 

\subsection{Design to achieve privacy}

Following this, we detail the process by which each non-faulty agent adds noise to its state in initial steps to ensure its initial state remains unrecoverable. The accumulated noise is then removed by the agent after those initial steps. 
Hence, this corresponds to changing~\eqref{eq:lti} to
\begin{equation}\label{eq:lti_aug_sel2}
\begin{array}{rcl}
\tilde{x}^{(k+1)}&=&\tilde{A}\tilde{x}^{(k)}+\mathbb I_{[0,T]}(k)\tilde{\eta}^{(k+1)}\\
{x}^{(k+1)}&=&P^{(k)}\tilde{x}^{(k+1)},
\end{array}
\end{equation}
where $\tilde{\eta}^{(k)}$ is the aforementioned noise added by the agents. For simplicity, the noise vector each entry of the noise vector $\tilde{\eta}^{(k)}$ is ${\eta}^{(k)}$, but it is easy to allow different noise values, ensuring the same theoretical guarantees. 
Moreover, 
\[\tilde{x}^{(k)} = \left[\begin{smallmatrix}x_{1,\mathcal S_1}^{(k)}&\ldots&x_{n,\mathcal S_1}^{(k)}&\ldots&x_{n,\mathcal S_1}^{(k)}&\ldots,x_{n,\mathcal S_r}^{(k)}\end{smallmatrix}\right]^\intercal.
\] 
Our proposed method is structured so that each agent keeps track of all possible ``partial'' consensus protocols considering the states' updates according to the possible subsets formed by parts of agents, which is captured in the augmented state $\tilde{x}^{(k)}$. The goal is to latter  differentiate/determine the subset of agents that are misbehaving 
(e.g., $\tilde{x}_{v,\emptyset}^{(k)}$ is the state where $\emptyset$ nodes are discarded, and $\tilde{x}_{v,\{1,3\}}^{(k)}$ is the state where nodes 1 and 3 are discarded), starting as~$x_v^{(0)}$. 

These steps are encapsulated in Algorithm~\ref{alg:gn}, in Section~\ref{sub:ifwtni}.

Next, we specify how each non-faulty agent  selects its current state from a vector of possible states, effectively excluding the states of faulty agents. In other words, we are defining $P^{(k)}$. The steps are mirrored in Algorithm~\ref{alg:ss}, in Section~\ref{sub:sss}. 


\subsection{State augmentation and initialization }\label{sub:sai}

Our proposed method is structured so that each agent keeps track of all possible ``partial'' consensus protocols considering the states' updates according to the possible subsets formed by parts of agents, which is captured in the augmented state $\tilde{x}^{(k)}$. The goal is to latter  differentiate/determine the subset of agents that are misbehaving 
(e.g., $\tilde{x}_{v,\emptyset}^{(k)}$ is the state where $\emptyset$ nodes are discarded, and $\tilde{x}_{v,\{1,3\}}^{(k)}$ is the state where nodes 1 and 3 are discarded), starting as~$x_v^{(0)}$. 

    Another important property we want to achieve is average consensus (not just consensus). Hence, a building block of our method is Algorithm~\ref{alg:ds} that allows agents to design the weights to average the received information from neighbors~\cite{bu2018accelerated}. 
    If the iteration number is smaller than $T$, a non-faulty agent randomly generates noise with a prescribed standard deviation $\xi$ and adds that value to the vector state.

\begin{algorithm}[!ht]
\caption{\textsc{Doubly\_Stochastic}}
\scriptsize
\begin{algorithmic}[1]
\STATE \textbf{Input:} $Adj$, $discard$
{\tiny\color{gray}  			 \hfill$\rhd$ $Adj$ is the adjacency matrix of the network (without self-loops) and $discard$ is a set of nodes to discard}
\STATE $A \gets -Adj$; $m \gets$ number of rows of $A$
\STATE $A[i][j] \gets 0$, $\forall_{(i\in discard)\lor(j\in discard)}$ {\tiny\color{gray}  			 \hfill$\rhd$ the discarded nodes have weights set to 0 so consider the subnetwork without them}
\STATE $d_{\max} \gets 1$
\FOR{$i$ \textbf{from} $1$ \textbf{to} $m$}
    \STATE $tmp \gets -\displaystyle\sum_{j=1}^{m}A[i][j]+1$
    \STATE $A[i][i] \gets tmp$
    \STATE $d_{\max} \gets (d_{\max} \text{ \textbf{if} } d_{\max} > tmp \text{ \textbf{else} } tmp$)
\ENDFOR
\FOR{$i$ \textbf{from} $1$ \textbf{to} $m$}
    \FOR{$j$ \textbf{from} $1$ \textbf{to} $m$}
        \STATE $A[i][j] \gets -\frac{A[i][j]}{d_{\max}}+\begin{cases}
        \frac{d_{\max}+1}{d_{\max}}, & \text{if }i=j\\
        0, & \text{otherwise}
        \end{cases}
        $
    \ENDFOR
\ENDFOR
\FOR{$i$ \textbf{in} $discard$}
    \STATE $A[i][i] \gets 0$
\ENDFOR
\RETURN $A$
\end{algorithmic}
\label{alg:ds}
\end{algorithm}

\begin{remark}\label{rem:doubly}
    Another distributed way of computing the agents weights that leads to a doubly-stochastic dynamic matrix~\cite{garin2010survey} can be used in the design of our proposal, i.e., replacing Algorithm~\ref{alg:ds}.~\hfill$\blacktriangle$
\end{remark}

\subsection{Initial finite window of time noise injection}\label{sub:ifwtni}
    
Next, the last building block of the method we propose is Algorithm~\ref{alg:gn}. This algorithm simply generates noise and adds the noise to the agent states, allowing each agent to keep track of the accumulated introduced noise. The accumulated noise serves to correct the consensus from past noise, allowing to reach average consensus.  
Since the noise is added by the different agents to their respective states, it is replicated across the entire augmented state, kept by the different agents. 

\begin{algorithm}[!ht]
\caption{\textsc{Gaussian\_Noise}}
\scriptsize
\begin{algorithmic}[1]
\STATE \textbf{Input:} $x$, $ps$, $k$, $\xi$, $i$, $ns$, $T$
{\tiny\color{gray}  			 \hfill$\rhd$ $x$ is an agent state vector at a  iteration $k$, $ps$ is $\mathcal P([n],0,f)$, $\xi$ is the noise standard deviation, $i$ is the agent, $ns$ is the accumulated noise, $T$ is the number of noisy iterations}
\IF{$k < T$}
    \STATE $\eta[k] \gets \text{normal}(0,\xi)$ {\tiny\color{gray}  			 \hfill$\rhd$ noise added up to a certain number of iterations}
\ELSE
    \STATE $\eta \gets -ns[i]$ {\tiny\color{gray}  			 \hfill$\rhd$ first iteration without noise cancels the accumulated noise, after that this value is zero}
\ENDIF
\STATE $ns[i] \gets ns[i] + \eta[k]$ {\tiny\color{gray}  			 \hfill$\rhd$ accumulated noise}
\FOR{$el$ \textbf{in} $ps$}
    \STATE \textit{last}($x[el]$) $\gets$ \textit{last}($x[el]$) $+ \eta[k]$ {\tiny\color{gray}  			 \hfill$\rhd$ add noise to the state \\[.1cm] \hfill vector the \textbf{\textit{last}} operator returns the last element of a list and $x[el]$\\[.1cm]
    \hfill  collects the states' evolution of the agent when excluding the\\
    \hfill  agents in $el$}
\ENDFOR
\end{algorithmic}
\label{alg:gn}
\end{algorithm}


    

\subsection{Scalar state selection}\label{sub:sss}

The intuition behind the proposed process is that entries in the augmented state vector of a non-faulty agent, which reflect the consensus of a subnetwork including faulty agents, asymptotically approach a convex combination of the asymptotic states of the faulty agents. Therefore, the algorithm identifies the first entry free from the influence of faulty agents, if such an entry exists. If no such entry is found, the algorithm defaults to using the entry that incorporates information from all agents. 

Specifically, the non-faulty agent $v$ select the entry of its augmented state to be its state value as the entry it is confident to be the one excluding faulty agents, see Lemmas~\ref{lemma:1} and~\ref{lemma:2}.  
This step is summarized in Algorithm~\ref{alg:ss}. 

In Algorithm~\ref{alg:ss}, we store the augmented states of each iteration of agent $v$ in $w$. Specifically, $w$ is a \textit{hashmap} (or \textit{dictionary}) data structure, where the keys are sets representing the excluded agents, and the values are lists containing the consensus values for each iteration. For example, $w[\mathcal{U}]$ holds the consensus evolution of agent $v$ when the agents in $\mathcal{U}$ are excluded.

\begin{algorithm}[!ht]
\caption{\textsc{Select\_State}}
\scriptsize
\begin{algorithmic}[1]
\STATE \textbf{Input:} $w$, $v$, $\mathcal S$, $\varepsilon$ 
{\tiny\color{gray}  			 \hfill$\rhd$ $w$ is a \textit{hashmap} with keys being the elements of $\mathcal S$ (see the details in the preamble of the pseudo-code), $v$ is an agent, and $\varepsilon$ is the precision used to compare values}
\FOR{$i$ \textbf{from} $1$ \textbf{to} $f$}
    \STATE $ps_i \gets \{x\,:\,  x\in \mathcal P(\mathcal S,i)\}$ {\tiny\color{gray}  			 \hfill$\rhd$ parts of size $i$}
    \STATE $\ell \gets \left|\{x\,:\, x\in \mathcal P(\mathcal S,i-1)\}\right|$ {\tiny\color{gray} \\ 			 \hfill$\rhd$ number of parts of size $< i$}
    \IF{$i=1 \bigwedge \displaystyle\mathop{\forall}_{\substack{s\in ps_i
    \\v\notin s
     }
    }\left|\text{\textit{last}}(w[ps_i[s]]\text{)}-\text{\textit{last}}(w[\emptyset]\text{)}\right| > \varepsilon$}
        \RETURN $\text{\textit{last}}(w[\emptyset]\text{)}$ {\tiny\color{gray}  			 \hfill$\rhd$ all the subnetworks lead to different consensus and thus there is no faulty agents detected}
    \ENDIF
    \FOR{$j$ \textbf{from} $1$ \textbf{to} $|ps_i|$}
        \IF{$
        \displaystyle\mathop{\forall}_{\substack{s=k\ldots,\ell+j\\ s\neq \ell+j}}\left|\text{\textit{last}}(w[ps_i[s]]\text{)}-\text{\textit{last}}(w[ps_i[j]]\text{)}\right| > \varepsilon
        $}
            \RETURN $\text{\textit{last}}(w[ps_i[j]]\text{)}$
            {\tiny\color{gray}  			 \hfill$\rhd$ the first entry that \\discards exactly the faulty nodes, the remaining ones have at least one faulty agent and converge to that agent's value}
        \ENDIF
    \ENDFOR
\ENDFOR
\RETURN $\text{\textit{last}}(w[\emptyset]\text{)}$ {\tiny\color{gray}  			 \hfill$\rhd$ no faulty agents detected}
\end{algorithmic}
\label{alg:ss}
\end{algorithm}

\subsection{Main algorithm} 

A non-faulty agent $v$ updates its vector state and repeats steps in Sections~\ref{sub:ifwtni} and~\ref{sub:sss} until reaching a stopping criteria. Using the state vectors received from its neighbors, the non-faulty agent $v$ updates each entry of its vector state using the respective entry from the neighbors and excluding the information from agents in the key that indexes each entry of its state vector. 

For example, $\tilde{x}_{v,\mathcal X}^{(k+1)}$ is updated using the values in $\left\{\tilde{x}_{u,\mathcal X}^{(k)}\,:\,u\in (\mathcal N_v\setminus\mathcal X)\cup\{v\}\right\}$. Note that an agent following this update rule is non-faulty and, therefore, it considers its current state information. 
These steps are detailed in Algorithm~\ref{alg:main}. Notice that, as we previously mentioned, the pseudo-code is not presented in a distributed setting for reproducibility reasons. Notwithstanding, the main for-loops are effortlessly implementable is a distributed fashion.



\begin{algorithm}[!ht]
\caption{\textsc{Average\_Consensus}}
\scriptsize
\begin{algorithmic}[1]
\STATE \textbf{Input:} $x_0$, $Adj$, $\xi$, $N$, $f$, $\mathcal F$, $\mathcal F\_vals$, $\varepsilon$, $T$
{\tiny\color{gray}  			 \hfill$\rhd$ $x_0$ is the initial agents' state, $Adj$ is the adjacency matrix of the agents' network, $\xi$ is the standard deviation of the noise that will be added, $N$ is the number of iterations, $f$ is the maximum number of permitted faulty agents, $\mathcal F$ is the set of faulty agents, $\mathcal F\_vals$ has the functions that dictate the faulty agents evolution, $varepsilon$ is the tolerance used to compare numbers, and $T$ is the number of noisy iterations $T<<N$}
\STATE $n \gets |x_0|$ {\tiny\color{gray}  			 \hfill$\rhd$ number of agents}
\STATE $ns \gets \underbrace{\{0,\ldots,0\}}_{n\text{ times}}$ {\tiny\color{gray}  			 \hfill$\rhd$ each agent accumulated noise}
\STATE $ps \gets \mathcal P([n], 0, f)$ {\tiny\color{gray}  			 \hfill$\rhd$ parts of agents of size up to $f$}
\STATE $A[s] \gets $ \textsc{Doubly\_Stochastic}$(Adj, s)$,  $\forall_{s\in ps}$ {\tiny\color{gray}  			 \hfill$\rhd$ weight matrix for each subnetwork excluding up to $f$ agents}
\STATE $X[i][s] \gets \left\{x_0[i]\right\}$ $\forall_{s\in ps} \forall_{i\in\left([n]\setminus\mathcal F\right)}$ {\tiny\color{gray}  			 \hfill$\rhd$ agents' initial vector state}
\FOR{$j$ \textbf{from} $1$ \textbf{to} $n$}
    \STATE \textsc{Gaussian\_Noise}($X[j]$, $ps$, $0$, $\xi$, $j$, $ns$, $T$) {\tiny\color{gray}  			 \hfill$\rhd$ add noise to each agent initial vector state}
\ENDFOR
\IF{$\mathcal F \neq \emptyset$}
   
            \STATE $X[j][s] \gets \{\mathcal F\_vals[j](0)\}$, $\forall_{s\in ps} \forall_{j\in\mathcal F}$
        \STATE $x[0][j] \gets \mathcal F\_vals[j](0)$, $\forall_{j\in\mathcal F}$
         {\tiny\color{gray}  			 \hfill$\rhd$ faulty agents set their initial vector state $X[j][s]$ and selected state $x[0][j]$}
\ENDIF
\FOR{$i$ \textbf{from} $1$ \textbf{to} $N-1$}
    \FOR{$j\in\left([n]\setminus\mathcal F\right)$}
            \FOR{$s$ \textbf{in} $ps$}
                \STATE $X[j][s][i] \gets \displaystyle\sum_{v\in\mathcal (N_j\setminus s)\cup\{j\}}X[v][s][i-1].A[s][j][v]$  {\tiny\color{gray}  			 \hfill$\rhd$ non-faulty agents follow the update rule}
            \ENDFOR
            \STATE \textsc{Gaussian\_Noise}($X[j]$, $ps$, $i+1$, $\xi$, $j$, $ns$, $T$) {\tiny\color{gray}  			 \hfill$\rhd$ non-faulty agents add noise to the vector state}
            \STATE $x[i][j] \gets  \textsc{Select\_State}(X[j][:][i], j, ps, f, \varepsilon)$ {\tiny\color{gray}  			 \hfill$\rhd$ non-faulty agents follow select the state as one of the vector state entries; \quad $v[:]$ denote all entries of $v$}
        \ENDFOR
        \FOR{$j\in\mathcal F$}
            \FOR{$s\in ps$}
                \STATE $X[j][s][i] \gets \mathcal F\_vals[j](i)$
            \ENDFOR
            \STATE $x[i][j] \gets \mathcal F\_vals[j](i)$ {\tiny\color{gray}  			 \hfill$\rhd$ faulty agents update their vector state and selected state}
        \ENDFOR
    \ENDFOR
\RETURN $x$ {\tiny\color{gray}  			 \hfill$\rhd$ all agents states for the $N$ iterations}
\end{algorithmic}
\label{alg:main}
\end{algorithm}

The presented pseudocode of our implementation can be easily made distributed. Our goal here is to make the reproducibility of the presented results a simple task\footnote{A GitHub repository with a Python 3 implementation can be found here: \href{https://github.com/xuizy/Average_consensus_with_resilience_and_privacy_guarantees_without_losing_accuracy}{GitHub code}.}. 

\begin{remark}\label{rem:directed}
    Observe that the proposed method works for directed networks with the proviso that all the subnetworks that need to be considered allow to design weights that yield a doubly-stochastic dynamic matrices.~\hfill$\blacktriangle$
\end{remark}

Next, we present several results that culminate in the soundness of Algorithm~\ref{alg:main}. The first result states that if an agent is considering (even if indirectly) information from faulty agents, then it converges to a convex combination of the faulty agents asymptotic state. In other words, the agent fails to achieve the desired consensus value. 

\begin{lemma}\label{lemma:1}
Consider Algorithm~\ref{alg:main}, the vector state of agent $v$, $\tilde{x}_v$ and let $\mathcal S\in\mathcal P([n], 0, f)$. If $|\mathcal S|<|\mathcal F|$ then 
\[
\displaystyle\lim_{k\to\infty}\tilde{x}_{v,\mathcal S}^{(k)} = \beta,
\]
where $\beta = \displaystyle\sum_{u\in\mathcal F}\alpha_u x_u^{(0)}$, $\alpha_u\geq 0$ and $\displaystyle\sum_{u\in\mathcal F}\alpha_u = 1$.~\hfill$\diamond$
\end{lemma}

\begin{myproof}
Consider the vector state of agent $ v $, $ \tilde{x}_v\in\mathbb R^{|\mathcal P([n],0,f)|} $, and let $ \mathcal{S} \in \mathcal{P}(\text{[n]}, 0, f) $. If $ |\mathcal{S}| < |\mathcal{F}| $, then there is at least one faulty agent in $ \mathcal{V} \setminus \mathcal{S} $.

By property (ii) in Assumption~\ref{ass:faulty}, the faulty agents asymptotically behave as stubborn agents. 
In this scenario, we can model the consensus evolution of the subnetwork $ \mathcal{H} = (\mathcal{V} \setminus \mathcal{S}, \mathcal{E}') $, a block of $\tilde A$, as

\begin{equation}\label{eq:lti2}
\left[\begin{smallmatrix}
x^{(k+1)}_R \\
x^{(k+1)}_S
\end{smallmatrix}\right]
= 
\left[\begin{smallmatrix}
Q_1 & Q_2 \\
0 & I
\end{smallmatrix}\right]   
\left[\begin{smallmatrix}
x^{(k)}_R \\
x^{(k)}_S
\end{smallmatrix}\right],
\end{equation}
for $ k \geq T $ (when no noise is added). This equation represents the system up to a label permutation of the agents, which corresponds to a permutation of rows and columns in the dynamic matrix and the state vector. 

The subvector $ x^{(k)}_R $ corresponds to the states of the non-faulty agents, and $ x^{(k)}_S $ corresponds to the states of the faulty (stubborn) agents. Note that these subvectors  correspond to subvectors of the vector $\tilde x$, which comprises the agents' augmented vector states. 
Moreover, a stubborn agent does not change its state and, therefore, the dynamics of stubborn agents is captured by the identity matrix. 

It is easy to see that $ Q_1 $ is a substochastic matrix since the columns sum up to a value less or equal to 1, with at least a row that sums up to less than 1. Furthermore, its spectral radius is $\rho(Q_1)<1$~\cite{horn2012matrix}, yielding an exponentially stable matrix~\cite{horn2012matrix}.





Note that $ x^{(T)}_S = x^{(0)}_S $. 
Furthermore, it is easy to prove by induction that

\begin{equation}\label{eq:stub4}
\left[\begin{smallmatrix}
Q_1 & Q_2 \\
0 & I
\end{smallmatrix}\right]^k  
= 
\left[\begin{smallmatrix}
Q_1^k & \quad & \left(\sum_{i=0}^k Q_1^i\right)Q_2 \\
0 & \quad & I
\end{smallmatrix}\right].
\end{equation}

So, we can pose \eqref{eq:stub4} as

\begin{equation}\label{eq:stub5}
\begin{cases}
x_R^{(k+1)} = Q_1^k x_R^{(T)} + \left(\sum_{i=0}^k Q_1^i\right)Q_2 x_S^{(k)}, \\[.2cm]
x_S^{(k+1)} = I x_S^{(k)}.
\end{cases}
\end{equation}

Next, since $ Q_1 $ is asymptotically stable, it follows that

\[
\lim_{k \to \infty} Q_1^k = \mathbf{0} \quad \text{and} \quad \lim_{k \to \infty} \left(\sum_{i=0}^k Q_1^i\right)Q_2 = (I - Q_1)^{-1} Q_2.
\]
We have that 
\begin{equation}\label{eq:stub6}
\begin{cases}
\displaystyle\lim_{k\to\infty} x_R^{(k)} = \lim_{k\to\infty}Q_1^k x_R^{(T)} + \lim_{k\to\infty}\left(\sum_{i=0}^k Q_1^i\right)Q_2 x_S^{(k)}, \\[.2cm]
\displaystyle\lim_{k\to\infty} x_S^{(k)} = \lim_{k\to\infty}I x_S^{(k)},
\end{cases}
\end{equation}
which is the same as
\begin{equation}\label{eq:stub7}
\begin{cases}
\displaystyle x_R^{(\infty)} =   (I - Q_1)^{-1} Q_2 x_S^{(\infty)}, \\[.2cm]
\displaystyle x_S^{(\infty)} = x_S^{(\infty)}.
\end{cases}
\end{equation}
Note that for a faulty agent $v$, its corresponding entry in $x_S^{(\infty)}$ is $\varphi_v$ (recalling Assumption~\ref{ass:faulty}). 
Hence, we conclude that the non-faulty agents' limit state is a convex combination of the faulty (stubborn) agents' initial states.
\end{myproof}

On the opposite, the second result states that if an agent is not considering (even if indirectly) information from faulty agents, then it converges to the average of the initial agents' states of the agents being considered. That is, the agent successfully achieves the desired consensus value.

\begin{lemma}\label{lemma:2}
    Consider Algorithm~\ref{alg:main}, the vector state of agent $v$, $\tilde{x}_v$ and let $\mathcal S\in\mathcal P([n], 0, f)$. If $\mathcal F\subset \mathcal S$ and $\mathcal H=(\mathcal V\setminus\mathcal S,\mathcal E')$ is connected, then 
\[
\displaystyle\lim_{k\to\infty}\tilde{x}_{v,\mathcal S}^{(k)} =\frac{1}{n-|\mathcal S|}\sum_{j\in \mathcal V\setminus\mathcal S} x_j^{(0)}.
\]
~\hfill$\diamond$
\end{lemma}

\begin{myproof}
    Consider an agent $v\in\mathcal V\setminus\mathcal S$. 
    For the first $T$ as in~\eqref{eq:lti_aug_sel2} steps where the agents add noise, let $A\equiv A[\mathcal S]$ built in Algorithm~\ref{alg:main}, then the dynamics presented in the algorithm regarding the entry $\tilde{x}_{v,\mathcal S}^{(k)}$  can be posed as 
    \begin{equation}\label{eq:noise_no_attack}
        \begin{cases}
            \displaystyle \tilde{x}_v^{(k+1)} =  \sum_{u\in(\mathcal N_v\setminus\mathcal S)\cup\{v\}}A^{\mathcal S}_{uv} \left(x_v^{(k)} + \eta_v^{(k+1)}\right), \\[.2cm]
            \displaystyle ns_v^{(k+1)} = ns_v^{(k)} + \eta_v^{(k+1)}.
        \end{cases}
    \end{equation}
    Because $A^{\mathcal S}$ is doubly-stochastic, we have that $\mathbf{1}^\intercal A^{\mathcal S} = \mathbf{1}^\intercal$, and since  $x^{(k+1)}=A^{\mathcal S}\left(x^{(k)}+\eta^{(k+1)}\right)$, it follows that 
    \begin{equation}\label{eq:um_transposto}
    \mathbf{1}^\intercal x^{(k+1)}=\mathbf{1}^\intercal\left(x^{(k)}+\eta^{(k+1)}\right)=\mathbf{1}^\intercal x^{(k)}+\mathbf{1}^\intercal \eta^{(k+1)}.
    \end{equation}
    From~\eqref{eq:um_transposto}, we have that 
    \begin{equation}\label{eq:um_transposto2}
    \mathbf{1}^\intercal x^{(k+1)}=\mathbf{1}^\intercal x^{(0)}+\mathbf{1}^\intercal\left(\sum_{i=0}^{k}\eta^{(k+1)}\right).
    \end{equation}
    Next, when $k = T$, with $ T = \left\lceil \frac{N}{\text{noisy\_prop}} \right\rceil $, $\eta^{(T)}=-\sum_{i=0}^{T-1}\eta^{(i)}$. Hence, using~\eqref{eq:um_transposto2}, we have that 
    \[\begin{split}
     \mathbf{1}^\intercal x^{(T+1)}  =  \mathbf{1}^\intercal x^{(T)}+\mathbf{1}^\intercal \eta^{(T)}\\
      = \left(\mathbf{1}^\intercal x^{(0)}+\mathbf{1}^\intercal\left(\sum_{i=0}^{T-1}\eta^{(i)}\right)\right)+ \mathbf{1}^\intercal\left(-\sum_{i=0}^{T-1}\eta^{(i)}\right)\\
      =  \mathbf{1}^\intercal x^{(0)}.
    \end{split}\]
    Therefore, at time $T+1$ the average of the  agents' states is the average of the initial agents' states. Hence, the agents converge to the average of the initial values. 
\end{myproof}

\begin{theorem}[Soundness]\label{th:sound}
    Under the Assumptions~\ref{ass:-1}--\ref{ass:faulty}, Algorithm~\ref{alg:main} reaches average consensus of non-faulty agents whenever $f\geq |\mathcal F|$.~\hfill$\circ$
\end{theorem}

\begin{myproof}
    To show that Algorithm~\ref{alg:main} reaches average consensus of non-faulty agents whenever $f\geq |\mathcal F|$, we first note that, by invoking Lemma~\ref{lemma:1}, we know that each subnetwork of the original network containing at least one agent in $\mathcal F$ converges to a convex combination of the faulty agents asymptotic states. 

    Moreover, by invoking Lemma~\ref{lemma:2}, we  know that each subnetwork of the original network not containing agents in $\mathcal F$ converges to the average of that subnetwork agents' initial states. 

    Therefore, it only remains to show that the state selection captured in Algorithm~\ref{alg:ss} effectively pinpoints the non-faulty agents state vector entry that corresponds to the subnetwork $\mathcal H=(\mathcal V\setminus\mathcal F,\mathcal E')$, i.e., the exact network that excludes the faulty agents and only those. 

    First, we note that Algorithm~\ref{alg:ss} selects as the agent state the first entry of the vector state that satisfy the predicate \begin{equation}\label{eq:th_main}     \displaystyle\mathop{\forall}_{\substack{s=k\ldots,k+j\\ s\neq k+j}}\left|\text{\textit{last}(w}[ps[s]]\text{)}-\text{\textit{last}(w}[ps_i[j]]\text{)}\right| > \varepsilon,
    \end{equation}
    where $k =\left|\{x\,:\, |x| < i, x\in \mathcal P(\mathcal V,0,f)\}\right|$.  If none of the entries satisfy the predicate, it selects the first entry, corresponding to the consensus state of the network with all the agents. 
    Note that $\varepsilon$ is arbitrarily small and corresponds to the precision utilized in the comparison. 
    In practice, if the computations are performed with numbers represented by a data type with a specified precision (e.g., 64-bit floating point numbers), then we may even think of $\varepsilon$ as $0$ with that precision. 
    
    So, if there are $|\mathcal F|=0$ faulty agents, we know that Algorithm~\ref{alg:ss} terminates in step~\textsf{A}\ref{alg:ss}.6 because the condition 
    \[  \displaystyle\mathop{\forall}_{\substack{s\in \mathcal P(\mathcal V,i)\\v\notin s}}\left|\text{\textit{last}(w}[s]\text{)}-\text{\textit{last}(w}[\emptyset]\text{)}\right| > \varepsilon.
    \]
    holds true when $N\to\infty$ by invoking  Lemma~\ref{lemma:2} and Assumption~\ref{ass:1}. 

    If $0<|\mathcal F|\leq f$, then we have that from~\eqref{eq:th_main}, when $N\to\infty$ the first the entry of a non-faulty agent that corresponds to the set that excludes exactly the agents in $\mathcal F$ is selected. 
    This property holds since any of the previous entries of that vector state have at least one attacked agent, and by Lemma~\ref{lemma:1} the value converges to a convex combination asymptotic faulty agents'  state. 
    
    Additionally, any other set of excluded agents $\mathcal S$ such that $|\mathcal S|= |\mathcal F|$ and $\mathcal S\neq \mathcal F$ leads to excluding strictly less than $|\mathcal F|$ faulty agents. Moreover, for each such set $\mathcal S$ there is another set $\mathcal S'$ such that $\mathcal S\cap\mathcal F = \mathcal S'\cap\mathcal F$ (by excluding distinct non-faulty agents) and yielding the save convex combination of faulty agents, recall Assumption~\ref{ass:faulty}. Thus, by invoking Assumption~\ref{ass:1}, it follows that Algorithm~\ref{alg:ss} terminates in step~\textsf{A}\ref{alg:ss}.10 with the correct state resulting from the average consensus of non-faulty agents. 
\end{myproof}

\begin{remark}
    It is worth noting that a value of $f\geq |\mathcal F|$ allows to achieve the desired resilient consensus. However, the larger the value of $f$, the more costly becomes the use of Algorithm~\ref{alg:main} (see Proposition~\ref{prop:complex}).~\hfill$\blacktriangle$
\end{remark}

{

Subsequently, we show that Algorithm~\ref{alg:main} successfully achieves privacy of the initial agents' states. We split the proof into two parts: (i) when no resilience is considered, which will be addressed in Theorem~\ref{th:private}, and (ii) when resilience is considered, which will be presented as a corollary in Corollary~\ref{th:privateAll}.

\begin{theorem}[Privacy without resilience]\label{th:private}
Consider the (undirected) network of agents $\mathcal G(A) = (\mathcal V, \mathcal E)$ satisfying Assumption~\ref{ass:priv} and the iteration scheme:
\begin{equation}
  x^{(k+1)} = Ax^{(k)} + I_{[0,1]}(k)\eta^{(k+1)},
\end{equation}
where $\mathcal G$ is connected and $A$ is doubly stochastic. 
Without loss of generality, each agent $i$ introduces private noise $\eta^{(0)}_i$ at $k=0$ and, for $T=1$, removes it at $k=1$ by setting $\eta^{(1)}_i = -\eta^{(0)}_i$. For any finite $T>1$, the same privacy guarantees hold by considering $\sum_{k=0}^{T-1}\eta^{(k)}_i$ as the initial noise and removing it at time $T$.
Let $x^{(0)} = \tilde{x}^{(0)} + \eta^{(0)}$ denote the noise-corrupted initial state that agents share with their neighbors, where $\tilde{x}^{(0)}$ is the private initial state. 
Then agent $u$ cannot uniquely determine the private initial state $\tilde{x}^{(0)}_v$ for any agent $v$ that either:
\begin{enumerate}
   \item[$(i)$] is not a neighbor of $u$ (i.e., $v \notin \mathcal{N}_u$), or
   \item[$(ii)$] is a neighbor of $u$ but has at least one neighbor outside $\mathcal{N}_u \cup \{u\}$ (as guaranteed by Assumption~\ref{ass:priv}).
\end{enumerate}
\end{theorem}

\begin{myproof}
    The proof establishes that under the given conditions, the noise-corrupted initial state $\tilde{x}^{(0)}_v$ cannot be uniquely determined from agent $u$'s measurements. We first establish the system structure and then demonstrate why privacy is preserved through unobservability analysis.

Consider the augmented system that explicitly captures both states and noise injection:
\begin{equation}
  \left[\begin{smallmatrix}
  \tilde{x}^{(k+1)}\\
  \eta^{(k+1)}
  \end{smallmatrix}\right] = 
  \left[\begin{smallmatrix} 
  A & I_n\\
  0_n & 0_n
  \end{smallmatrix}\right]
  \left[\begin{smallmatrix}
  \tilde{x}^{(k)}\\
  \eta^{(k)}
  \end{smallmatrix}\right] = \bar{A}\bar{x}^{(k)},
\end{equation}
where $\tilde{x}^{(0)} = x^{(0)} + \eta^{(0)}$ represents the noise-corrupted initial states that are shared among neighbors. The noise is removed after one step by setting $\eta^{(1)} = -\eta^{(0)}$, ensuring convergence to the true average while maintaining privacy of initial conditions.

We have that agent $u$ receives measurements described by
\begin{equation}
  y^{(k)} = \left[\begin{smallmatrix}C_1 & 0 \\ 0 & C_2 \end{smallmatrix}\right]
  \left[\begin{smallmatrix}
  \tilde{x}^{(k)}\\
  \eta^{(k)}
  \end{smallmatrix}\right] = \bar{C}\bar{x}^{(k)},
\end{equation}
where $C_1 \in \mathbb{R}^{(|\mathcal{N}_u|+1) \times n}$ captures state observations of agent $u$ and its neighbors $\mathcal{N}_u$, while $C_2 \in \mathbb{R}^{1 \times n}$ captures noise observation only for agent $u$.

For this system, any least squares estimate $\hat{\bar{x}}^{(0)}$ of the initial augmented state can be written as
\begin{equation}
  \hat{\bar{x}}^{(0)} = \bar{x}^{(0)} + U_wv,
\end{equation}
for some vector $v$, where $U_w$ spans the unobservable subspace of $(\bar{A},\bar{C})$, \textit{i.e.}, the unobservable subspace corresponding to agent $u$'s measurements. 

The corresponding estimate of the private initial state is:
\begin{equation}
  \hat{x}^{(0)} = H\bar{x}^{(0)} + HU_wv,
\end{equation}
where $H = [I_n \;\; I_n]$ maps the augmented states to private states.

The privacy of agent $v$ is violated if and only if $e_v^\top HU_w = 0$. This follows from the least squares estimation framework: any estimate of the initial augmented state has the form $\hat{\bar{x}}^{(0)} = \bar{x}^{(0)} + U_wv$ for arbitrary $v$, where $U_w$ spans the unobservable subspace. Then $\hat{x}^{(0)}_v = e_v^\top H\bar{x}^{(0)} + e_v^\top HU_wv$, making unique determination of the initial state possible if and only if $e_v^\top HU_w = 0$.

Hence, to characterize $U_w$, we analyze the PBH matrix pencil equation. The unobservable subspace $U_w$ consists of vectors $\bar{x} =[x^\top \; \eta^\top]^\top$ that satisfy:
\begin{equation}\label{eq:nullspace_augmented}
\left[\begin{smallmatrix}
\lambda I_n - A & -I_n\\
0_n & \lambda I_n \\
C_1 & 0 \\
0 & C_2
\end{smallmatrix}\right]
\left[\begin{smallmatrix}
x\\
\eta
\end{smallmatrix}\right] = 
0,
\end{equation}
for some value of $\lambda$.
In fact, for $\lambda \neq 0$, the second block row requires $\lambda\eta = 0$, forcing $\eta = 0$. This case corresponds to dynamics without privacy-inducing noise, so we focus on $\lambda = 0$, i.e., when noise affects the initial condition.
This is the non-trivial case when the noise-free consensus protocol is observable from agent $u$'s measurements.

For $\lambda = 0$, we get
\[
  Ax + \eta = 0\quad\land\quad
  C_1x = 0\quad\land\quad
  C_2\eta = 0.
\]
Subsequently, substituting $\eta = -Ax$, we obtain the equivalent condition 
\begin{equation}\label{eq:nullspace_reduced}
  \left[\begin{smallmatrix}
  C_1\\
  C_2A
  \end{smallmatrix}\right]x = 
  \left[\begin{smallmatrix}
  C_1\\
  e_w^\top A
  \end{smallmatrix}\right]x = 0.
\end{equation}
Therefore, privacy is preserved if and only if $e_v^\top H\bar{x} = 0$ for any $\bar{x}$ that solves~\eqref{eq:nullspace_augmented}. Equivalently, the latter condition can be reformulated to $e_v^\top(I-A)x = 0$ for any vector $x$ satisfying~\eqref{eq:nullspace_reduced}. This condition is equivalent to having $e_v^\top (I-A) $ being linearly independent from the rows of $\left[\begin{smallmatrix} C_1 \\ e_w^\top A \end{smallmatrix}\right]$. This linear independence means there exists no linear combination of these rows that could reveal $v$'s noise-corrupted initial state. 

To establish this linear independence, we analyze the structure of these row vectors:
\begin{itemize}
   \item [$1)$] $e_v^\top(I-A)$ is row $v$ of the graph Laplacian, with non-zero elements at positions corresponding to $v$ and its neighbors.

\item[$2)$] $e_w^\top A$ is row $w$ of $A$, and is linearly dependent on $C_1$ rows since both only involve observations from $\mathcal{N}_u \cup \{u\}$.

\item[$3)$] The rows of $C_1$ are exactly $e_j^\top$ for $j \in \mathcal{N}_u \cup \{u\}$.
\end{itemize} 
Finally, the linear independence of $e_v^\top(I-A)$ from $C_1$ rows follows from the following two cases:
\begin{itemize}
   \item [$1.$]For $v \notin \mathcal{N}_u$: $e_v^\top(I-A)$ contains non-zero elements at positions of $v$ and its neighbors, some of which are not observed by $u$.
   \item[$2.$] For $v \in \mathcal{N}_u$, by Assumption~\ref{ass:priv}, $v$ has a neighbor $j \notin \mathcal{N}_u \cup \{u\}$, creating a non-zero element in $e_v^\top A$ at position $j$ that cannot be generated by $C_1$ rows.
\end{itemize}
This linear independence ensures $e_v^\top HU_w \neq 0$, proving that agent $u$ cannot uniquely determine the noise-corrupted initial state $\tilde{x}^{(0)}_v$ for any agent $v$ satisfying either condition of the theorem.
\end{myproof}

}

The result in Theorem~\ref{th:private} entails the following corollary. 

\begin{corollary}[Privacy with resilience]\label{th:privateAll}
Given the (undirected) network of agents $\mathcal G = (\mathcal V, \mathcal E)$,  consider the iteration scheme of~\eqref{eq:lti}, where the dynamics matrix $A$ is primitive and doubly stochastic. Assume that a set of nodes $\mathcal P$ follows the predefined update rule. Then, for any resilient parameter $f$, if the noise each agent adds is kept private, a
malicious/curious agent $u$ cannot recover the initial value of another agent  $x_v^{(0)}$, for $v\in \mathcal P$, there is $w\in\mathcal V$ with $u\neq v$ such that $u\notin\mathcal N_v$ in any subnetwork of $\mathcal G$ that excludes at most $f$ agents.  
\end{corollary}

\begin{myproof}
By Assumption~\ref{ass:faulty}, each subnetwork of agents that consists of excluding up to $f$ agents remains strongly connected and, by Assumption~\ref{ass:priv}, we can invoke Theorem~\ref{th:private} to each of these subnetworks, entailing the desired result. 
\end{myproof}

\begin{remark}\label{rem:5}
Note that adding noise in an initial finite window of time may lead to a temporary erroneous classification of a non-faulty agent as being faulty. However, after the initial window of noise, the classification is sound with the proviso that the number of faulty-agents does not exceed the parameter $f$.~\hfill$\blacktriangle$
\end{remark}

Remark~\ref{rem:5} enables privacy and resilience without trade-offs, unlike differential privacy approaches. 
We now analyze Algorithm~\ref{alg:main}'s computational complexity. 

\begin{proposition}[Computational complexity]\label{prop:complex}
    A non-faulty agent $v$ following the average consensus update rules of Algorithm~\ref{alg:main}, in a network with $n$ agents, with neighbors $\mathcal N_v$, maximum number of allowed faulty agents $f$, and number of iterations $N$  incurs in a computational time-complexity of $\mathcal O(\max\{|\mathcal N_v|, f n^{2f}\}N).$~\hfill$\bullet$
\end{proposition}

\begin{myproof}
   To compute the computational complexity that non-faulty agent $v$ incurs, we assess the cost of each step in Algorithm~\ref{alg:main} from the viewpoint of a single agent.  
   Steps~\textsf{A}\ref{alg:main}.1--\textsf{A}\ref{alg:main}.2 have a constant cost. 
   step~\textsf{A}\ref{alg:main}.3 has $\mathcal O(n^f)$ cost and step~\textsf{A}\ref{alg:main}.4, from agent $v$ viewpoint, only needs only to check the neighbors of agent $v$ to compute the weights of $A$ related to $v$ in Algorithm~\ref{alg:ds}, and repeat it for $\mathcal O(n^f)$ elements, yielding a cost of $\mathcal O(n^f|\mathcal N_v|)$. 
   
   Next, in step~\textsf{A}\ref{alg:main}.8, the agent generates random noise and adds it to its vector state, with Algorithm~\ref{alg:gn}, with cost $\mathcal O(n^f)$. 
   Subsequently, in step~\textsf{A}\ref{alg:main}.17, the agent computes a weighted average of the received states with cost $\mathcal O(|\mathcal N_v|)$, and it adds noise in step~\textsf{A}\ref{alg:main}.19 with cost $\mathcal O(n^f)$.  
   These steps (steps~\textsf{A}\ref{alg:main}.17--\textsf{A}\ref{alg:main}.19) are repeated $\mathcal O(n^f)$ times, yielding a cost of $\mathcal O(\max\{|\mathcal N_v|,n^{2f}\})$. 
   In step~\textsf{A}\ref{alg:main}.20, it selects an element of its vector state to be its state, using Algorithm~\ref{alg:ss}. 
   
   The cost upper bound is $\mathcal O( (n^f)^2 f)=\mathcal O(f n^{2f})$. 
   Last, steps~\textsf{A}\ref{alg:main}.17--\textsf{A}\ref{alg:main}.20 are repeated $\mathcal O(N)$ times. 
   In summary, the total cost is the maximum of the previous costs, which is $\mathcal O(\max\{|\mathcal N_v|, f n^{2f}\}N)$.
\end{myproof}

\begin{remark}
    Observe that the computational time-complexity an agent incurs in when using Algorithm~\ref{alg:main} update rules and disables resilience ($f=0$) is $\mathcal O(|\mathcal N_v| N)$, which is the cost of the usual consensus method.~\hfill$\blacktriangle$
\end{remark}

{
\section{Illustrative Examples}\label{sec:examples}

We illustrate the use of the proposed average consensus method (Algorithm~\ref{alg:main}) with some examples. 
For visualization purposes, in the subsequent examples, we replace step~\textsf{A}\ref{alg:gn}.3 of Algorithm~\ref{alg:gn} with the following two steps
\[
    \begin{split}
    \theta \gets & \text{normal}(0,\xi)\\
    \eta \gets & \begin{cases}
         \theta, & \text{if }ns[i] \leq 0\\
         -\theta, & \text{otherwise.}
     \end{cases}
     \end{split}
\]
By making this small adjustment, the added noise will not become overpowering, preserving the readability of the illustrative examples in the plots.  

\begin{figure}[!ht]
    \centering
 \includegraphics[width=.33\textwidth]{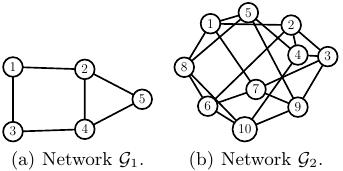}
    \vspace{-1mm}
    \caption{Networks of agents utilized in Examples I--II.}
    \label{fig:grafos}
\end{figure}

In the following section, we present illustrative examples of consensus networks. 

\subsection{Example I - a faulty agent in a 5 agents network with non-constant state (ensuring privacy)}

The following example uses the same setup as the previous one, with $x_0=\left[\begin{smallmatrix}
    f_1(0) & 0.4 & 0.65 & 0.1 & 0.8
\end{smallmatrix}\right]$. The faulty agents are $\mathcal F=\{1\}$ and $x_1(k)=f_1(k)=\frac{0.6 \cos (k)}{k+1}$. 
In this case, the average consensus must be $0.4875$, as achieved by the agents running Algorithm~\ref{alg:main} with $f=1$, see Fig.~\ref{fig:examplo_3}. In this case, the number of iterations is $N=90$, and there are $T=10$ noisy iterations.

\begin{figure}[!ht]
    \centering
\includegraphics[width=.9\textwidth]{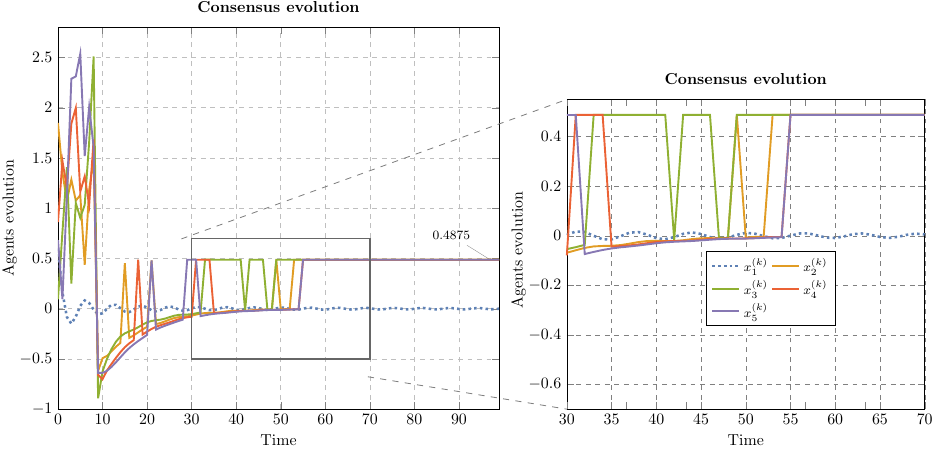}
\vspace{-1mm}
    \caption{Agents consensus evolution of Example I, with network of agents $\mathcal G_1$, and $\mathcal F=\{1\}$.}
    \label{fig:examplo_3}
\end{figure}

Moreover, the average consensus \textbf{ensures privacy}. For network $\mathcal G_1=(\mathcal V,\mathcal E)$, no agent has all others as neighbors: $\mathcal N_1=\{2,3\}$, $\mathcal N_2=\{1,4,5\}$, $\mathcal N_3=\{1,4\}$, $\mathcal N_4=\{2,3,5\}$, and $\mathcal N_5=\{2,4\}$. The subnetworks excluding one agent (Fig.~\ref{fig:subnetsG2}) maintain this property. Thus, by Theorem~\ref{th:privateAll}, agents' initial states remain private.

\begin{figure}[!ht]
    \centering
  \includegraphics[width=.33\textwidth]{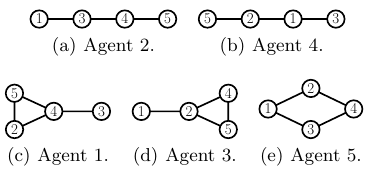}
    \vspace{-1mm}
    \caption{Subnetworks of $\mathcal G_1$ excluding one agent (indicated in each subfigure).}
    \label{fig:subnetsG2}
\end{figure}




\subsection{Example II - 3 faulty agents in a 10 agents network sharing asymptotically distinct values (ensuring privacy)}

Finally, to illustrate Algorithm~\ref{alg:main} with a larger network and more faulty agents, we utilized the network $\mathcal G_2$ with set of faulty agents $\mathcal F=\{1,2,7\}$.
In this example, the agents initial states are $$x_0=\left[\begin{smallmatrix}
    f_1(0) & f_2(0) & f_3(0) & 0.606 & 1.295 & 1.753 & 0.0807 & 0.593 & 1.657 & 0.549
\end{smallmatrix}\right]^\intercal,$$ where $x_1(k)=f_1(k)=0.2+2\frac{\sin{x}}{x+1}$, $x_2(k)=f_2(k)=0.6-3\frac{\cos{x}}{x+1}$, and $x_3(k)=f_3(k)=0.8-3\frac{\cos{x}}{x+1}$ ($\lim_{k\to\infty}f_1(k)=0.2$, $\lim_{k\to\infty}f_2(k)=0.6$ and $\lim_{k\to\infty}f_3(k)=0.8$). 
Thus, the desired consensus is $\approx 0.9333$, as the non-faulty agents successfully converge to, when using Algorithm~\ref{alg:main} with $f=3$, see Fig.~\ref{fig:examplo_5}. In this case, the number of iterations is $N=80$, and there are $T=10$ noisy iterations.

\begin{figure}[!ht]
    \centering
\includegraphics[width=.9\textwidth]{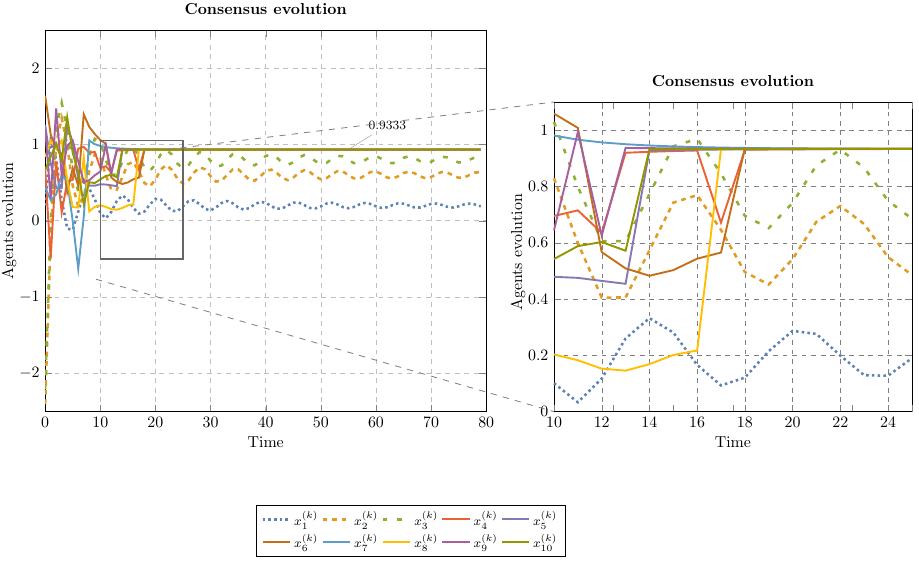}
    \caption{Agents consensus evolution of Example VI, with network of agents $\mathcal G_2$, and $\mathcal F=\{1,2,3\}$.}
    \label{fig:examplo_5}
\end{figure}

To verify that the initial states of the agents remain private, consider the network $\mathcal G_2 = (\mathcal V, \mathcal E)$, where each agent $v \in \mathcal V$ has exactly $|\mathcal N_v| = 4$ neighbors. Consequently, Theorem~\ref{th:private} applies, ensuring privacy for each agent. Furthermore, Theorem~\ref{th:privateAll} also holds for any subnetwork of $\mathcal G_2$ that excludes up to three agents. This is straightforward to verify, as the subnetworks range from 7 to 10 agents, meaning that excluding up to three agents will never result in a subnetwork where a single agent has all the remaining agents as neighbors. Thus, \textbf{privacy is preserved} even in these smaller subnetworks. 

Note that, in this example, there are 176 subnetworks, and we cannot do the analysis cases-by-case as we did in the previous example.

}











\section{Conclusions}\label{sec:conclusion}

In this paper, we address private and resilient average consensus in discrete-time networked systems. While traditional methods require trade-offs between accuracy, privacy, and resilience, our proposed consensus mechanism integrates these factors without compromise, albeit at increased computational and communication costs.
Recent work achieved similar results but required calculating the left-eigenvector of the system dynamics matrix, necessitating full knowledge of agents' weights or centralized coordination. Our approach eliminates these requirements. By introducing controlled noise into agents' states, we mask individual contributions while maintaining malicious activity detection and accurate average consensus.
The framework's polynomial time complexity relative to network size and attack resilience demonstrates its practicality, supported by empirical evidence and numerical results.



{\small
\bibliography{bib.bib}
}

\end{document}